\documentclass{amsart}

\usepackage{amsmath, amsthm, amsopn, amsfonts}

\usepackage{color}

\usepackage{graphicx}

\usepackage{graphpap}

\def\pasdegrille{\let\grille = \pasgrille}

\def\aat#1#2#3{

\divide \dimen1 by 48

\dimen3=\dimen1

\multiply \dimen1 by #1

\advance \dimen1 by -\dimen3

\divide \dimen1 by 101

\multiply \dimen1 by 100

\divide \dimen2 by \count11

\multiply \dimen2 by #2

\setbox0=\hbox{#3}\ht0=0pt\dp0=0pt

  \rlap{\kern\dimen1 \vbox to0pt{\kern-\dimen2\box0\vss}}\dimen1= \wd1

\dimen2=\ht1}

\def\pasgrille{

\count12= \dimen1

\divide \count12 by 50

\divide \dimen2 by \count12

\count11 =\dimen2

\

\divide \dimen1 by 48

\setlength{\unitlength}{\dimen1}

\smash{\rlap{\ }}

\dimen1= \wd1

\dimen2=\ht1

}

\def\grille{

\count12= \dimen1

\divide \count12 by 50

\divide \dimen2 by \count12

\count11 =\dimen2

\

\divide \dimen1 by 48

\setlength{\unitlength}{\dimen1}

\smash{\rlap{\graphpaper[1](0,0)(50, \count11)}}

\dimen1= \wd1

\dimen2=\ht1

}

\pasdegrille

\usepackage{amssymb}

\ifx\pdfoutput\undefined

\usepackage{graphicx}  \DeclareGraphicsExtensions{.ps}

\else

\usepackage{amsfonts}

\usepackage[pdftex]{graphicx}  \DeclareGraphicsExtensions{.pdf,.mps}

\fi

\usepackage[T1]{fontenc}

\usepackage{epic, eepic}

\def\squarebox#1{\hbox to #1{\hfill\vbox to #1{\vfill}}}

\newcommand{\1}{{\bold 1}}

\newcommand{\CI}{{\mathcal C}^\infty }

\newcommand{\CIc}{{\mathcal C}^\infty_{\rm{c}} }

\newcommand{\SP}{{\mathbb S}}

\newcommand{\NN}{{\mathbb N}}

\newcommand{\rest}{\!\!\restriction}

\theoremstyle{plain}

\newtheorem{thm}{Theorem}

\newtheorem{prop}{Proposition}[section]

\theoremstyle{definition}

\newtheorem{rem}{Remark}

\numberwithin{equation}{section}

\title[Eigenfunctions for partially rectangular billiards]
{Eigenfunctions for partially rectangular billiards}
\author[J. Marzuola]{Jeremy Marzuola}
\address{Mathematics Department, University of California \\
Evans Hall, Berkeley, CA 94720, USA}
\email{marzuola@math.berkeley.edu}

\usepackage{amssymb}
\usepackage{amsmath, amsthm, amsopn, amsfonts}

\def\11{{\rm 1~\hspace{-1.4ex}l} }
\def\R{\mathbb R}

\begin{document}

\maketitle   
   
\section{Introduction}   
\label{in}
In this note, we further develop the methods of Burq-Zworski \cite{BZ3} 
to study eigenfunctions for billiards which have rectangular components: these include the Bunimovich 
billiard, the Sinai billiard, and the recently popular pseudointegrable billiards
\cite{Bog}.  The results are an application of a "black box" point of view as
presented in \cite{BZ2} by the same authors.

\begin{figure}[ht]
\includegraphics[width=10cm]{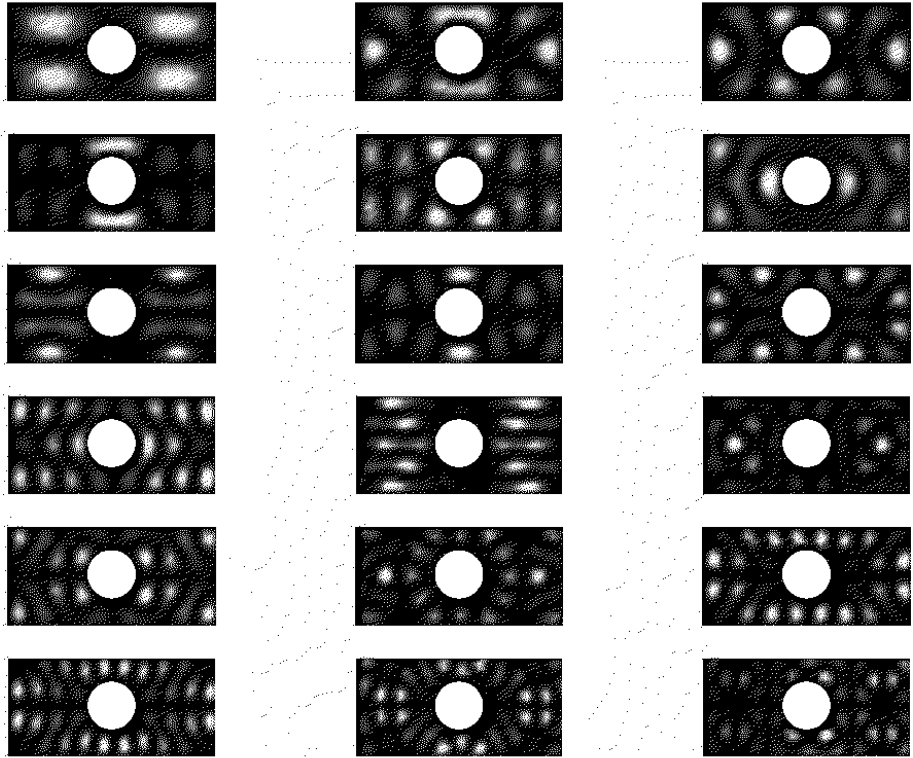}
\caption{Experimental images of eigenfunctions in a Sinai billiard
microwave cavity -- see {\tt http://sagar.physics.neu.edu}. We see
that there is always a non-vanishing presence near the boundary of the
obstacle as predicted by Theorem \ref{t:s} below.}
\label{fig:sin}
\end{figure}

By a partially rectangular billiard, we mean a connected planar domain,
$ \Omega $, with a piecewise smooth boundary, which contains a  rectangle,
$ R \subset \Omega $, such that if we decompose the boundary of 
$ R $, into pairs of parallel segments, $ \partial R = \Gamma_1 \cup \Gamma_2 $,
then $ \Gamma_i \subset \partial \Omega $, for at least one $ i$.

We show that for such billiards, 
the eigenfunctions of the Dirichlet, Neumann,
or periodic Laplacian cannot concentrate in closed sets in the interior of the rectangular part.
A combination of this elementary result with 
the now standard, but highly non-elementary, 
propagation results of Melrose-Sj\"ostrand \cite{MeSj} and Bardos-Lebeau-Rauch \cite{BLR},
can give further improvements -- see \cite{BZ2},\cite{BZ3}.

Here, we prove further non-concentration results, away from the obstacle in the Sinai billiard
(see Fig.1 and Theorem \ref{t:s}), and along certain trajectories
in pseudointegrable billiards, (see Fig.5 and Theorem \ref{t:j}).  
For recent motivation coming from the study of {\em quantum chaos}
we suggest \cite{Bog},\cite{BZ3},\cite{Do},\cite{Ze}, and references given there.

\medskip
\noindent
{\sc Acknowledgments.} This paper is a development of an unpublished work by N. Burq and M. Zworski, 
who treated the case of a square billiard rather than a general rectangular billiard.
The author is very grateful to Maciej Zworski and Nicolas Burq for allowing me to use results from their 
unpublished work, as well as many helpful conversations.  I would also like to thank Srinivas Sridhar for 
allowing the use the experimental images shown in Fig.\ref{fig:sin}.

\section{Semiclassical Pseudodifferential Operators on a Torus}
\label{PDO}

In this section, we would like to discuss properties of Pseudodifferential Operators (PDO's) on a torus.  
To begin, we examine the nature of PDO's and their symbol classes.  In ${\mathbb R}^n$, we define a
Weyl quantization of an operator $a(x,hD)$ where $a \in {\mathcal S} ({\mathbb R}^{2n}), a=a(x, \xi)$ by:
$$a (x,hD) u(x) = \frac{1}{(2 \pi h)^n} \int_{{\mathbb R}^n} \int_{{\mathbb R}^n} e^{\frac{i}{h} \langle x-y, \xi \rangle}
a(x, \xi) u(y) dy d\xi ,$$
for $u \in {\mathcal S}$ and a symbol class by:  
$$S^k_\delta (m) = \{ a \in C^{\infty} ({\mathbb R}^{2n}) \| | \partial^\alpha a | \leq C_\alpha h^{-\delta |\alpha| -k} m 
\ \text{for all multi-indices} \ \alpha \},$$
where $m: {\mathbb R}^{2n} \to (0,\infty)$ is is an order function, i.e. there exist constants $C$, $N$ such that
$$ m(z) \leq C \langle z - w \rangle^N m(w).$$
We also have $$S^{-\infty}_\delta (m) = \bigcap_{k=-\infty}^{\infty} S^k_\delta (m).$$ 
For $k,\delta = 0$ we write simply $S (m)$. 
On a torus, however, $a$ and all its derivatives are bounded in the $x$ variable, thus for $h$ small and $k$ positive, we need not worry
about the derivatives in x, only those in $\xi$.  Also, for $k$ negative, provided
that we have the proper local regularity for our symbol $a$, this definition still works perfectly on a torus.

Note also that we need only work with symbols that are periodic in the x-variable
with period determined by the dimensions of the torus.  In other words,
$a(x,\xi)=a(x+\gamma, \xi)$ for $\gamma \in (a{\mathbb Z}) \times (b{\mathbb Z})$, where $a, b \in {\mathbb R}$.  
Using this relation, we easily see the following propostion.

\begin{prop}
\label{p:3}  
If $a(x,\xi)$ is a periodic symbol in x with period $\gamma$, then $a(x,D) T_\gamma = T_\gamma a(x,D)$ where
$T_\gamma u(x) = u(x - \gamma)$.  
\end{prop}
\begin{proof}
We easily calculate:

\begin{eqnarray}
\label{eq:p3}
a(x,h D) T_\gamma (u) & = & \int_{{\mathbb R}^n} a(x,\xi) e^{\frac{i}{h} <x-y,\xi>} \hat{u} (y-\gamma) dy d\xi \\
&  = & \int_{{\mathbb R}^n} a(x-\gamma,\xi) e^{\frac{i}{h} <x-\gamma-{\tilde y},\xi>} \hat{u} (\tilde y) d{\tilde y} d\xi \\
&  = & T_\gamma(a(x,D) u)
\end{eqnarray}
\end{proof}

From the above proposition, it becomes clear that the properties of symbol classes in Euclidean space  
translate directly to properties of similarly defined symbol classes on a torus.  
For instance, we have the following result.

\begin{prop}
\label{p:4}
Given $u(x) = u(x+\gamma)$ where $\gamma$ is as above, and $u$ is $L^2$ on a torus, then $a(x,h D) u(x)$ is
$L^2$ on the torus, when $a \in {\mathcal S}_\delta (1)$, $0 \leq \delta \leq \frac{1}{2}$.
\end{prop}
\begin{proof}
Note that the condition on $a$ implies that it is $L^2$ bounded.
Given a function $u(x)$ which is periodic on a torus, we can write it as 
$\sum_\gamma T_\gamma u_0$ where $u_0 = \chi (x) u(x)$ and $\chi (x)$ is equal to 1 on a single copy
of the torus in the plane and 0 otherwise.  Note that no assumptions about the smoothness of $\chi(x)$ are made.  
Hence, $u_0 \in L^2$ and therefore, so is $a(x,h D) u_0$.  Then, 
$a(x,h D) u(x) = \sum_\gamma T_\gamma a(x,h D) u_0$.  The sum converges for each $x$ since $a(x,D) u_0$ will
have compact support and we have
$a(x,h D) u$ a periodic function that is $L^2$ on the torus. 
\end{proof}

Using similar techniques, we would like to develop the concept of a microlocal defect measure in this setting.
As shown in \cite{Zw}, we have the following theorem in Euclidean space:
\begin{thm}
\label{t:dm}
There exists a Radon measure $\mu$ on ${\mathbb R}^n$ and a sequence $h_j \rightarrow 0$ such that 
\begin{equation}
\langle a^w (x,h_j D) u(h_j), u(h_j) \rangle \rightarrow \int_{{\mathbb R}^{2n}} a(x,\xi) d\mu
\end{equation}
for all symbols $a \in S(1)$.
\end{thm}

We call $\mu$ a microlocal defect measure associated with the family $\{ u(h) \}_{0<h \leq h_0}$. 
Note that an $S(1)$ symbol on the torus corresponds to an $S(1)$ symbol on the plane, therefore this result proves
the existence of microlocal defect measures on a torus as well.  

\begin{proof}
1.  Let $\{a_k \} \in C_c^\infty$ be dense in $C_0 ({\mathbb R}^{2n})$.  Select a sequence $h_j^1 \rightarrow 0$ such that
$$\langle a_1^w (x, h_j^1 D) u(h_j^1), u(h_j^1) \rangle \rightarrow \alpha_1.$$  Then, select a subsequence $\{ h_j^2 \} \subset \{ h_j^1 \}$
such that $$\langle a_2^w (x, h_j^2 D) u(h_j^2), u(h_j^2) \rangle \rightarrow \alpha_2.$$  Continue such that at the kth step, you take a subsequence
$\{ h_j^k \} \subset \{ h_j^{k-1} \}$ such that $$\langle a_k^w (x, h_j^k D) u(h_j^k), u(h_j^k) \rangle \rightarrow \alpha_k.$$  Then by
a diagonal argument, arrive at a sequence $h_j$ such that $$\langle a_k^w (x, h_j D) u(h_j), u(h_j) \rangle \rightarrow \alpha_k$$ for all $k=1,2,...$.

2.  Define $\Phi(a_k) = \alpha_k$.  By a standard theorem on operator norms, we have for each $k$ that 
$$| \Phi (a_k) | = | \alpha_k | = \lim_{h_j \rightarrow \infty} | <a_k^w u(h_j), u(h_j) > | \leq
\limsup_{h_j \rightarrow \infty} C \| a_k^w \|_{L^2 \rightarrow L^2} \leq C \sup |a_k|.$$  The mapping $\Phi$ is bounded, linear
and densely defined, therefore uniquely extends to a bounded linear functional on $S(1)$, with the estimate $$ | \Phi (a) | \leq C \sup |a|$$
for all $a \in S(1)$.  The Riesz Representation Theorem therefore implies the existence of a (possibly complex valued) Radon measure on
${\mathbb R}^{2n}$ such that $$\Phi (a) = \int_{{\mathbb R}^{2n}} a(x,\xi) d\mu.$$
\end{proof}

We now quote a general theorem about microlocal defect measures on Euclidean space which we can
then apply to a torus.  To state the propagation theorem in the form sufficient for our applications,
we follow \cite{BuIMRN}.

Let us consider a Riemannian manifold without boundary, $M$. By partitions of unity we can define semi-classical 
pseudo-differential operators $a(x, hD_x)$ associated to symbols $a(x, \xi)\in \CIc (T^* M)$. 
\par
Now we consider a sequence $(u_n)$ bounded in $L^2(M)$. satisfying
\begin{equation}\label{eq.equation}
(-h_n^2\Delta -1)u_n = 0.
\end{equation}  
Using~\eqref{eq.equation},
as in~\cite{GeLe93} (see also~\cite{BuIMRN}) we can prove the following result.

\begin{prop}
\label{prop2.1} There exist a subsequence $(n_{k})$ and a positive Radon  measure on $T^*M$, $\mu$ (a semi-classical measure for the sequence $(u_n)$), 
such that for any $a \in \CIc (T^* M)$
\begin{equation}
\label{eq2.8}\lim_{k\rightarrow + \infty}\left(a^w (x, h_{n_{k}}D_{x})u_{n_{k}}, u_{n_{k}}\right)_{ L^2(M)}= \langle \mu, a(x, \xi)\rangle.
\end{equation}
Furthermore this measure satisfies
\begin{enumerate}
\item The support of $\mu$ is included in the characteristic manifold:
\begin{equation}
\Sigma \stackrel{\rm{def}}{=} \{(x, \xi)\in T^*M; p(x,\xi)=\|\xi\|_x=1\}
\end{equation}
where  $\|\cdot \|_x$ is the norm for the metric at the point $x$,
\item The measure $\mu$ is invariant by the bicharacteristic flow (the flow of the Hamilton vector field of $p$):
\begin{equation}
\label{eq.inv}H_p \mu =0,
\end{equation}
\item For any $\varphi\in \CIc (T^*M)$, 
\begin{equation}
\label{eq.lim} \lim_{k \rightarrow + \infty} \|\varphi u_{n_k}\|^2 = \langle \mu, |\varphi|^2\rangle.
\end{equation}
\end{enumerate}
\end{prop}
The first two properties above are weak forms of the elliptic regularity and propagation of singularities results,  
whereas the last one states that there is no loss of $L^2$-mass at infinity in the $\xi$ variable.   

\begin{proof}
We will prove this proposition only for the case of a torus, but the methods are applicable to any manifold.
\bigskip

(0)  (Positivity)
We need to show that $a \geq 0$ implies $$\int_{{\mathbb T}^{2} \times {\mathbb R}^{2}} a(x, \xi) d\mu \geq 0.$$
Since $a \geq 0$, using the Garding inequality, we see that:  $$a^w (x,hD) \geq -Ch.$$
Let $h = h_j \rightarrow 0$, to see:
$$\int_{{\mathbb T}^{2} \times {\mathbb R}^{2}} a d\mu = \lim_{j \rightarrow \infty} \langle a^w (x,h_j D) u(h_j), u(h_j) \rangle \geq 0.$$

(1)  (Support of $\mu$)
Let $a$ be a smooth function such that $\text{supp} (a) \cap p^{-1} (1) = \emptyset$.  We must show 
$$\int_{{\mathbb T}^{2} \times {\mathbb R}^{2}} a d\mu = 0.$$
Select $\chi \in C^{\infty}_c ( {\mathbb T}^{2} \times {\mathbb R}^{2} )$ such that 
$$\text{supp} (a) \cap \text{supp} (\chi) = 0.$$
Then, $$a^w(x,hD) \left( ((p-1)^w+i \chi^w)^{-1} (p-1)^w \right) (x,hD) = a^w(x,hD) + {\mathcal O} (h^\infty)_{L^2 \rightarrow L^2}.$$
Apply $a^w (x,hD)$ to $u(h)$ to see that $a^w (x,hD) u(h) = o(1)$ and thus $\langle a^w (x,hD) u(h), u(h) \rangle \rightarrow 0$.
But, $$\langle a^w (x,hD) u(h_j), u(h_j) \rangle \rightarrow \int_{{\mathbb T}^{2} \times {\mathbb R}^{2}} a d\mu.$$

(2)  (Flow Invariance)
Select $a$ as above, then 
\begin{eqnarray} 
\langle [p^w,a^w] u(h), u(h) \rangle & = & \langle (p^w a^w -a^w p^w) u(h), u(h) \rangle \\
 & = & \langle a^w u(h), p^w u(h) \rangle - \langle p^w u(h), (a^w)^* u(h) \rangle \\
 & = & o(h) \ \text{as} \ h \rightarrow 0.
\end{eqnarray}
However, $[p^w,a^w] = \frac{h}{i} \{ p^w,a^w \} + O(h^2)$.  Hence, 
$$\langle [p^w,a^w] u(h), u(h) \rangle = \frac{h}{i} \langle \{ p,a \}^w u(h), u(h) \rangle + \langle o(h) u(h), u(h) \rangle.$$
As we let $h_j \rightarrow 0$, we get:  $$\int_{{\mathbb T}^{2} \times {\mathbb R}^{2}} \{p,a\} d\mu = 0.$$  
So, if $\Phi_t$ is the flow generated by the Hamiltonian vector field $H_p$, then 
$$\frac{d}{dt} \int_{{\mathbb T}^{2} \times {\mathbb R}^{2}} (\Phi_t* a) d\mu = 
\int_{{\mathbb T}^{2} \times {\mathbb R}^{2}} (H_p a)(\Phi_t) d\mu = 
\int_{{\mathbb T}^{2} \times {\mathbb R}^{2}} \{ p,a \} d\mu = 0.$$
Now, (3) follows easily by looking at the operator $|\varphi (x, \xi)|^2$ and applying the result about existence of a microlocal defect measure.
\end{proof}

\section{Partially rectangular billiards}
\label{prb}

In this section we will need to recall the basic control results \cite{Bu92},\cite{BZ2}
for rectagles, and the propagation results \cite{MeSj},\cite{BLR},\cite{BuIMRN},\cite{BuGe} for 
billiards. Since in the specific application presented in Section \ref{app} we only use propagation
away from the boundary, that is the only case we will review.

The following result from ~\cite{Bu92} is related to some earlier control results 
of Haraux~\cite{Ha} and  Jaffard~\cite{Ja}.\footnote{We remark that as noted in \cite{Bu92} 
the result holds for 
any product manifold $M= M_{x}\times M_{y}$, and the proof is essentially the same.}
\begin{prop}
\label{p:1} Let $\Delta$ be  the Dirichlet, Neumann, or periodic
Laplace operator on the rectangle $R= [0, 1]_{x} \times [0,a]_{y}$. 
Let $\omega_x$ be a non-empty open subset of $[0,1]$.  Then for any non-empty  $\omega \subset R$ of the form $
\omega= \omega_{x} \times [0,a]_{y}$, there 
exists $C$ such that for any solutions of
\begin{equation}
(\Delta -z) u =f \ \text{ on $R$}, \ u \rest_{\partial R}=0
\end{equation}
we have
\begin{equation}
\label{eq:6.12}\|u\|^2_{{L^2(R)}}\leq C \left(\|f \|^2_{H^{-1} (
[0,1]_{x}; L^2([0,a]_{y})) } +
\|u\rest _{\omega}\|^2_{{L^2(\omega)}} \right)
\end{equation}
\end{prop}
\begin{proof}
We will consider the Dirichlet case (the proof is the same in the other two 
cases) and 
decompose $u,f$ in terms of 
the basis of $L^2([0,a])$ formed by the Dirichlet eigenfunctions
$e_{k}(y)=  { \sqrt {{2}/a}}\sin(2k\pi y/a)$,
\begin{equation}
u(x,y)= \sum_{k}e_{k}(y) u_{k}(x), \qquad f(x,y)= \sum_{k}e_{k}(y) f_{k}(x).
\end{equation}
We get for $u_{k}, f_{k}$ the equation
\begin{equation}\label{estres.1}
\left(\Delta_{x}-\left(z+\left({2k\pi}/{a}\right)^2\right)\right)u_{k}= f_{k},\qquad u_{k}(0)=u_k(1)=0.
\end{equation}
We now claim that 
\begin{equation}
\label{eq:cont}
\|u_{k}\|^2_{{L^2([0,1]_{x})}}\leq C \left(\|f_k \|^2_{H^{-1}([0,1]_{x})} +
\|u_k \rest _{\omega_{x}}\|^2_{{L^2(\omega)}}\right), 
\end{equation}
from which, by summing the squares in $k$, we get~\eqref{eq:6.12}.

To see \eqref{eq:cont} we can use the propagation result above \ref{prop2.1} in dimension
one, but in this case an elementary calculation is easily available -- see
\cite{BZ3}.
\end{proof}

The following theorem is an easy consequence of Proposition \ref{p:1}:

\begin{thm}
\label{t:1}
Let $ \Omega $ be a partially rectangular billiard with the rectangular
part $ R \subset \Omega $, $ \partial R = \Gamma_1 \cup \Gamma_2 $, a decomposition 
into parallel components satisfying $ \Gamma_2 \subset \partial \Omega $. 
Let $ \Delta $ be the Dirichlet or Neumann Laplacian on $\Omega $. Then
for {\em any} neighbourhood of $ \Gamma_1 $ in $ \Omega $, $ V $, there exists $ C $
such that
\begin{equation}
\label{eq:t1}
- \Delta u = \lambda u \ \Longrightarrow \  \int_V | u ( x ) |^2 dx 
\geq \frac1C \int_R | u ( x ) |^2 dx \,, 
\end{equation}
that is, no eigenfuction can concentrate in $ R $ and away from $ \Gamma_1$.
\end{thm}

\begin{proof}
Let us take $x,y$ as 
the coordinates on the stadium, so that $x$ parametrizes $ \Gamma_2 \subset \partial 
\Omega $ and $ y$ parametrizes $ \Gamma_1 $, then  
\[ R = [0,1]_{x}\times [0,a]_{y} \,.\]
Let $\chi\in \CIc((0,1))$ be equal to $1$ on $[\varepsilon, 1- \varepsilon]$. Then $\chi(x) u(x,y)$ is a solution of
\begin{equation}
(\Delta-z)\chi u = [\Delta, \chi] u\text{ in $R$}
\end{equation} with the boundary conditions satisfied 
on $\partial R$. Applying Proposition~\ref{p:1}, we get 
\begin{equation}
\|\chi u\|_{L^2( R)}\leq C \left \| [ \Delta, \chi] u \|_{H^{-1}_{x}; L^2_{y}}+ 
\| u\rest _{\omega_{\varepsilon}}\|_{L^2( \omega_{\varepsilon})}\right)
\leq C' \| u\rest _{\omega_{\varepsilon}}\|_{L^2( \omega_{\varepsilon})} \,,
\end{equation}
where $\omega_{\varepsilon}$ is a neighbourhood of the support of $ \nabla\chi$.
Since a neighbourhood of $ \Gamma_1 $ in $ \Omega $ has to contain $ \omega_\varepsilon $
for some $ \varepsilon $, \eqref{eq:t1} follows.
\end{proof}

\section{Applications}
\label{app}

\begin{figure}
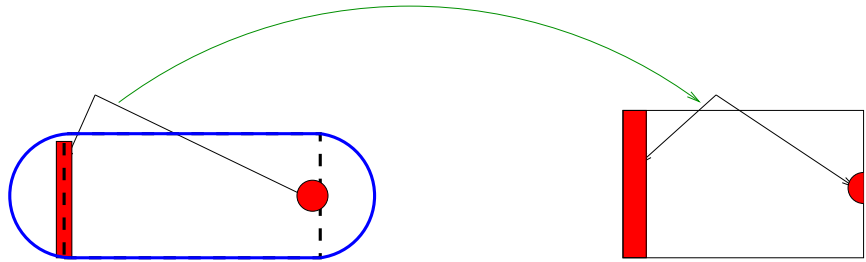

\include{figure2}
\caption{Control regions in which eigenfunctions have positive mass and the
rectangular part for the Bunimovich stadium.}
\label{fig:bath}
\end{figure}

In \cite{BZ2} and \cite{BZ3}, Proposition \ref{p:1} is used to prove that in the
case of the Bunimovich billiard shown in Fig.\ref{fig:bath},
the states have nonvanishing 
density near the vertical boundaries of the rectangle. That follows from Theorem
\ref{t:1} which shows that we must have positive density in the wings of the
billiard, and the propagation result (in the boundary case) based on the fact that
any diagonal controls a disc geometrically (see  \cite[Section 6.1]{BZ2}; in fact 
we can use other control regions as shown in Fig.\ref{fig:bath}). 
Here we consider another case which accidentally generalizes a control theory 
result of Jaffard \cite{Ja}.

The Sinai billiard (see Fig.\ref{fig:sin}) 
is defined by removing a strictly convex open set, $ \mathcal O$, with 
a $ \CI $ boundary,
 from a flat torus, $ {\mathbb T}^2_{a,b} \stackrel{\rm{def}}{=}
 (a \SP^1 ) \times (b \SP^1 )$:
\[ S  \stackrel{\rm{def}}{=} {\mathbb T}^2_{a,b} \setminus {\mathcal O} \,.\]
The following theorem results by applying Theorem 1 to a torus with sides of 
arbitrary length.

\begin{thm}
\label{t:s}
Let $ V $ be any open neighbourhood of the convex boundary, $ \partial {\mathcal O} $, in 
a Sinai billiard, $ S $. If $ \Delta $ is the Dirichlet or Neumann
Laplace operator on $ S $ then there exists a constant, $ C = C ( V ) $, such that
\begin{equation}
\label{eq:ts}
- h^2 \Delta u = E(h) u \ \Longrightarrow \  \int_V | u ( x ) |^2 dx 
\geq \frac1C \int_S | u ( x ) |^2 dx \,, 
\end{equation}
for any $h$ and $|E(h) - 1| < \frac{1}{2}$.
\end{thm}
\begin{proof}
First note that we can easily limit ourselves to the case where our flat torus has
one side of length 1 and one side of length a.  Suppose that the result is not true, 
in other words, there exists a sequence of 
eigenfunctions $ u_n $, $ \| u _n \| = 1 $, with the corresponding 
eigenvalues $ \lambda_n \rightarrow \infty $,  such that $ \int_V | u_n ( x )|^2 dx
\rightarrow 0 $. 

We first observe that the only directions in the support of 
the corresponding semi-classical defect measure, $\mu$, have to be "rational", 
in other words, the trajectory must travel along a line of slope $\frac{m a}{n}$ where $m,n \in \NN$.  
The projection of a trajectory with an irrational direction is dense on the 
torus and hence must encounter the obstacle $ \partial {\mathcal O} $ (and 
consequently $V$). The propagation result recalled in Proposition \ref{prop2.1}, part (3),
gives a contradiction by choosing a proper test function $\phi$ which is nonzero on the support
of the measure $\mu$ resulting from our sequence of eigenfunctions (remark that we apply this result as long as the 
trajectory does not encounter the obstacle and consequently we need only the 
{\em interior} propagation).

Hence let us assume that there exists a rational direction in the support of 
the measure which then contains the periodic trajectory in that direction.  
As shown in Fig. \ref{f:3} we can find a maximal rectangular neighbourhood of 
the projection of that trajectory which avoids the obstacle.

\begin{figure}
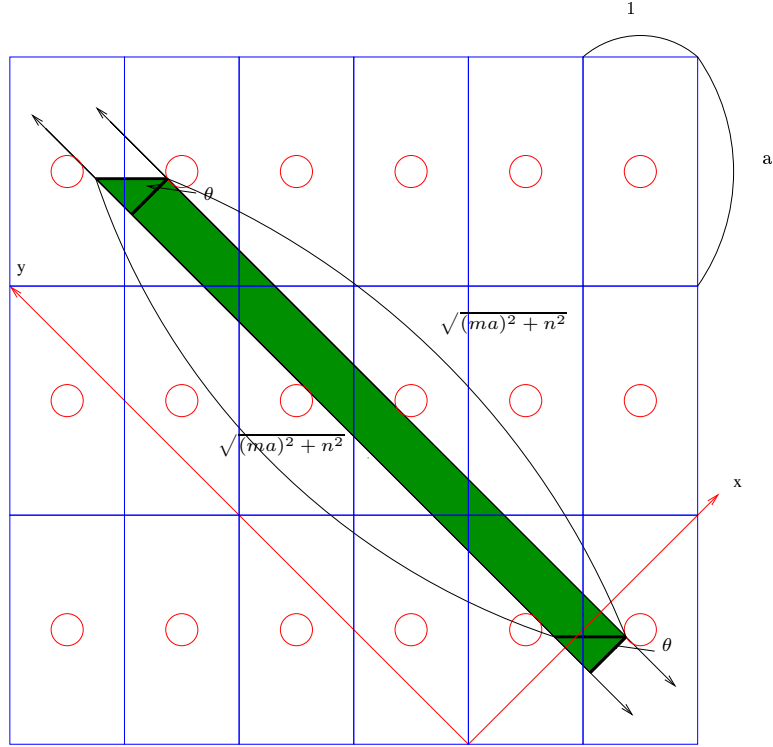

\include{figure3}
\caption{A maximal rectangle in a rational direction, avoiding the obstacle.
Because the parrallelogram is certainly periodic and our region has uniform width,
it is clear that the resulting rectangle is periodic.}
\label{f:3}
\end{figure} 

The rectangle can be described as $ R = [0,a_1]_{x_1} \times [0, b_1]_{y_1} $ with the the 
$ y_1 $ coordinate parametrizing the trajectory. 
Let $\epsilon$, $\delta >0$ be small.
Let $ u $ be an eigenfunction in our sequence and define 
$$ \chi = \{ \chi ( x_1 ) \in \CIc ({\mathbb T}^2) \ | \ \chi (x_1) = 1 \ \text{for all} \ x_1 \in (\epsilon,a_1-\epsilon), 
\ \chi(x_1) = 0 \ \text{outside $(0,a_1)$} \}.$$  
Note that we can then write $\chi = \chi(x,y)$ for $(x,y) \in {\mathbb T}^2_{1,a}$ as $x_1$ is simply a rotation and translation of 
the standard coordinates.
Then $ \chi ( x,y ) u ( x, y ) $ is a function on all of $ R $ satisfying the periodicity condition. 
Let  $\Phi_\xi ( \nu ) = \Phi(\xi - \nu)$, where we define
$$\Phi = \{ \Phi (\xi) \in \CIc ({\mathbb R}^2) \ | \ \Phi (\xi) = 1 \ \text{for} \ \xi \in B(0,\delta), 
\ \Phi (\xi) = 0 \ \text{for} \ \xi \in {\mathbb R}^2 \setminus B(0,2 \delta) \}, $$ 
where $B(0,\delta)$ is a ball centered at $0$ of radius $\delta$.
Note, due to the compact support of this function in $\xi$,
$\Phi_\xi (D)$ is in the symbol class $S(\langle \xi \rangle^{-N})$ for any $N$.
Let $ \Delta_R $ is the (periodic) Laplacian on $ R $. Using Fourier decomposition
we can arrange that $ [ \Delta_R , \Phi_\xi (D) ] = 0 $.  Since our eigenfunction is only defined on 
${\mathbb T}^2_{1,a} \setminus {\mathcal O}$, let us introduce a smooth function ${\chi}_0$ which is $0$ on a neighborhood $U$
of
the obstacle and $1$ on ${\mathbb T}^2_{1,a} \setminus V$ where $U \subset V$.
 Choose $U$ and $V$ such that $\chi \chi_0 = \chi$.
Hence, 
\[  ( -h^2 \Delta_R - E(h) ) \Phi_\xi(D) \chi {\chi}_0 u = [ -h^2 \Delta_R , \Phi_\xi(D) \chi] u 
=\Phi_\xi (D) [ -h^2 \Delta_R , \chi] \chi_0 u + {\mathcal O} ( h^{ \infty } ) 
\,, \ \  \| u \| = 1 \,, 
\]
by the properties of $S(\langle \xi \rangle^{-N})$ operators acting on $L^2$ functions on a torus.
As in the proof of Proposition \ref{p:1}, we now see that 
\begin{equation}
\label{eq:mic}
 \| \Phi_\xi \chi {\chi}_0 u \|_{L^2} \leq C \int_\omega | \chi_0 u |^2 
+ {\mathcal O} ( h^{\infty } ) \,,
\end{equation}
where $ \omega $ is a neighbourhood of $ \nabla \chi $ (in
the calculus of semi-classical pseudo-differential operators). 
Since the semi-classical defect measure of $ \Phi_\xi \chi {\chi}_0 u $ (which is $|\Phi_\xi \chi {\chi}_0|^2\times \mu$) was assumed to be
non-zero, \eqref{eq:mic} shows that the measure of $ {\chi}_0 u \rest_\omega $
is non zero and consequently there is a point in the intersection of the supports of 
$\mu$ and ${\chi}_0 u \rest_\omega$. But $\mu$ is invariant
by the flow (as long as it does not intersect the obstacle) and hence, once we choose $\epsilon$, $\delta$ small enough 
such that all the cut-offs above are very close to the
boundary of $ R$, its support can be made intersect any neighbourhood of $ \partial {\mathcal O} $.
\end{proof}

Now, from the above theorem, we see the following simple, but important consequences:

\begin{rem}
Let $S= {\mathbb T}^2_{a,b} \setminus {\mathcal O}$ where $ {\mathcal O}$ is sufficiently smooth in the case
of Neumann boundary conditions, but otherwise lacking restrictions.  Then, for $V$ any open neighborhood of 
$\partial {\mathcal O}$, and $u$ a solution of $-h^2 \Delta u = E(h) u$ as above, then (4.1) is satisfied.  
This follows from the above argument as the convexity of the obstacle was never used.
Thus, the result holds for any obstacle (even connectedness is not assumed here)
and is applicable to the special case of pseudointegrable billiards (see for instance
\cite{Bog} for motivation and description).  In the next section, we use an argument similar to that above 
in order to say even more about concentration along trajectories in specific pseudointegrable billiards. 
By an elementary reflection principle, the result also holds for an obstacle inside a square with Dirichlet 
or Neumann conditions on the boundary of the square.
\end{rem}

\begin{rem}
The proof above gives in fact the following estimate for any open neighbourhood, say $V$, of the obstacle:
\begin{equation}
\begin{gathered}
\exists C; \forall u,\, f \in L^2(S) \text{ solutions of } (-\Delta + \lambda)u =f, \qquad u \rest_{\partial S} =0\\
\|u\|_{L^2(S)}\leq C \left( \|f\|_{L^2(S)}+ \|u \11 _V\|_{L^2(V)}\right)
\end{gathered} 
\end{equation}
and according to~\cite[Theorem 4]{BZ2}, this implies that the Schr{\"o}dinger equation in $S$ 
is exactly controllable by $V$ in finite time. In fact, by working on the time evolution equation, 
we could strengthen this result allowing an arbitrarily small time. 
\end{rem}
This latter result was previously known~\cite{Ja} 
for the particular case $\Theta= \emptyset$ ($S= \mathbb{T}^2$) but the proof was 
based on subtle results about Fourier series \cite{Ka}.

\def\cprime{$'$}

\begin{rem}
As shown in \cite[Theorem 2$'$]{BZ2}, the results of Ikawa and G\'erard on 
scattering by two convex obstacles (see \cite{BZ2} and references given there)
give an estimate on the maximal concentration of an eigenfunction (or a 
quasimode) on a closed orbit in a Sinai billiard. Let $ \chi \in \CI(S;[0,1])  $  be 
supported in  a small neighbourhood of a closed transversally reflecting 
orbit. Then  for any family $ ( - \Delta - \lambda ) u_\lambda = {\mathcal O}
( \lambda^{-\infty } ) $, $ \| u_\lambda \| = 1 $, 
\[ C \int_S  | u ( x ) |^2 ( 1 - \chi ( x ) ) dx \geq \frac{1}{ \log \lambda } \,, \]
that is a concentration on a closed trajectory, if at all possible, has to 
be very weak.
\end{rem}

\section{Pseudointegrable Billiards}
\label{PIB}

\begin{figure}
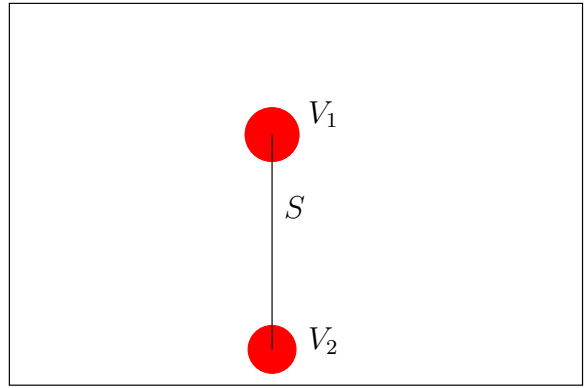

\include{figure4}
\caption{A pseudointegrable billiard $P$ consisting of a torus with a slit, $S$ along which we have Dirichlet boundary conditions.
We would like to show that eigenfunctions of the Laplacian on this torus must have concentration in the shaded regions $V_1$ and $V_2$.}
\label{f:4}
\end{figure}

We define a pseudointegrable billiard to be a plane polygonal billiard with corners whose angles are of the form $\frac{\pi}{n}$, 
for any integer $n$ (see \cite{Bog1}).  
In particular, we will be working with the billiard $P = {\mathbb T}^2_{a,b} \setminus { S}$ where
$ S$ is a slit that is parrallel to a side of the torus but not a closed loop.  In Remark 1, 
we point out that Theorem \ref{t:s} allows us to make statements about the $L^2$ mass
of eigenfunctions in a neighborhood of the slit for pseudointegrable billiards.  For this particular type of billiard,
it would be ideal to state that every eigenfunction must have non-zero mass in a small neighborhood of the edges of the slit (see Fig. \ref{f:4}).
In this section, we prove a weaker result about non-concentration along certain classical
trajectories in $P$ of semiclassical defect measures obtained from eigenfunctions $u$ such that $(-\lambda -\Delta)u = 0$ on $P$.

As with the Sinai billiard, the classical behavior of trajectories must be taken into account in our treatment of this problem.   
There cannot be concentration along trajectories that do not hit the slit as shown by Theorem \ref{t:s}.  If a trajectory 
has irrational slope, it is dense in $P$, and thus has mass near the edges of the slit as in Section \ref{app}.
Therefore, for our purpose, we concern ourselves only with rational trajectories which intersect the slit at some point.  As we are dealing with
periodic boundary conditions, let us consider the
plane tiled with copies of the billiard $P$.  

Assume that $P$ is oriented such that $S$ is parrallel to the $y$-axis.  
Let $\gamma \in S^* (P)$ be a trajectory.  Given the natural projection
$$ \pi_1: S^* (P) \rightarrow P,$$
we take $\gamma' = \pi_1 (\gamma)$, or the physical path along which the trajectory travels. 
Consider the projection
$$\tilde{\pi} : \R^2 \rightarrow {\mathbb T}^2_{a,b}.$$
We see that 
$$\tilde{\pi} : \R^2 \setminus \tilde{S} \rightarrow P, \ \text{where} \ \tilde{S} = \tilde{\pi}^{-1} (S).$$
Define $$\pi_2 : S^* (\R^2 \setminus \tilde{S} ) \rightarrow S^* (P)$$
to be the obvious projection.  Let $\tilde{\gamma} = \pi_2^{-1} (\gamma)$.
We can write $$\tilde{\gamma} = \displaystyle\cup_{j=1}^{\infty} \gamma_j,$$ where each $\gamma_j$ is a trajectory 
in $ S^* (\R^2 \setminus \tilde{S})$.  We note that by construction, $\gamma_i \cap \gamma_j = \emptyset$ for $ i \neq j$.  To see this, assume that 
$\gamma_i \cap \gamma_j = (x, \xi)$.  Then, $\gamma_i = \gamma_j$ as they would represent trajectories which travel 
through the same point in the same direction.  Now, let $$\pi^*_1 : S^* (\R^2 \setminus \tilde{S}) \rightarrow \R^2 \setminus \tilde{S}.$$
Select one trajectory from the above union, say $\gamma_1$.  Let $\gamma_1' = \pi^*_1 (\gamma_1)$. 
We see that either $\gamma_1'$ is bounded in the $x$-direction
or $\gamma_1'$ is unbounded in the $x$-direction.  Note that this property then holds for all $\gamma_j$, $j \in {\mathbb N}$.
For a trajectory $\gamma$, if the resulting path $\gamma_1'$ is bounded in the $x$-direction,
we say $\gamma$ is $x$-bounded.  We define $\gamma$ as $x$-unbounded if $\gamma_1'$ is unbounded in the $x$-direction.
See Fig. \ref{f:7} for examples.  Now, we are prepared to state our theorem concerning the billiard $P$.

\begin{figure}
\include{figure7}
\caption{Some examples of $x$-bounded and $x$-unbounded trajectories.}
\label{f:7}
\end{figure}

\begin{thm}
\label{t:j}
Let $\gamma$ be an $x$-bounded trajectory on $P= \mathbb{T}^2 \setminus S$.
If $ \Delta $ is the Dirichlet  Laplace operator on $ P $ then 
there exists no microlocal defect measure obtained from the eigenfunctions  on $P$ 
such that ${\rm{supp}} \  (d\mu) = \gamma$.
\end{thm}

\begin{proof}

Let  $\gamma '$  be as above.  Let $V_\epsilon$ be an $\epsilon$ neighborhood of $\gamma'$.  
If the theorem were false, we would have a sequence of eigenfunctions $u_n$, $\| u_n \|_{L^2} =1$ with the 
property
$$\displaystyle\int_{P \setminus V_\epsilon} | u_n |^2 dx \rightarrow 0,$$ for any $\epsilon.$ We show that
this is impossible.

For each $u_n$, we have $(-\Delta - \lambda_n) u_n = 0$, $u_n|_S =0$, $u_n \in L^2(P)$.  
Let $\tilde{\pi}$ be as above.
We define the sequence
$\tilde{u}_n = \tilde{\pi}^{-1} u_n$.  We have $(-\Delta - \lambda_n ) \tilde{u}_n = 0$,
$\tilde{u}_n|_{\tilde{S}} = 0$, and $\tilde{u}_n \in L^2_{\text{per}} (\R^2 \setminus \tilde{S})$.

If $\pi_2 : S^* (\R^2 \setminus \tilde{S}) \rightarrow S^* (P)$ is as above and 
$\tilde{\gamma} = \pi_2^{-1} (\gamma)$, then $\tilde{u}_n \rightarrow d\tilde{\mu}$ with 
$$\text{supp} \ (d\tilde{\mu}) = \tilde{\gamma} \subset S^* (\R^2 \setminus \tilde{S})\,. $$

Now, let $\pi_1^* : S^* (\R^2 \setminus \tilde{S}) \rightarrow \R^2 \setminus \tilde{S}$ be as above.
Select one trajectory, say $\gamma_1$.  As 
$\gamma_1$ is $x$-bounded, ${\gamma}_1' = \pi_1^* ( \gamma_1 )$ is contained in a strip in the plane which is 
infinite in the $y$-direction and bounded in the $x$-direction.
Thus, $\gamma_1'$ is contained in a strip, $C_0$, with minimal width in the $x$-direction.  Then,
$\tilde{u}_n$ satisfies $(-\Delta -\lambda_n) \tilde{u}_n = 0$ on $C_0$, is periodic in the $y$-direction, and
satisfies the following boundary conditions in the $x$-direction:  Dirichlet boundary conditions 
along the slits that intersect the boundary of $C_0$ and periodic boundary conditions otherwise.

Without loss of generality, we can choose the $x$-coordinates such that the 
boundaries of $C_0$ are
$x=-R$ and $x=0$.  We can then reflect to a strip, say $\tilde{C}_1$, with boundaries $x=-R$ and $x=R$, by defining a new function on $\tilde{C}_1$ by:
$$ \tilde{u}_n^{(1)} (x,y) =  \left\{ 
\begin{array}{ll}
 \ \ \tilde{u}_n (x,y)  &  x \in [-R,0], \\
 -\tilde{u}_n (-x,y)  & x \in (0,R). 
\end{array} \right.$$
Note that $\tilde{u}_n$ is periodic with period $2R$.
As a result, we have
$$(-\Delta - \lambda_n) \tilde{u}_n^{(1)} = f_n^{(1)}$$ on $\tilde{C}_1$, where
$$f_n^{(1)} = 2 u(0,y) {\delta}_{0} ' (x) - 2 u(R,y) {\delta}_{R} ' (x).$$
We note that $ f_n^{(1)} $ is supported away from the slits, $ \tilde S $.

Define 
$$\pi^{\sharp}_1 (x,y) = \left\{
 \begin{array}{ll}
 (x,y) & -R \leq x \leq 0, \\
  ( -x,y) & 0 \leq x \leq R.
  \end{array} \right.$$
If $ \pi_1^\sharp \; : \; \tilde{C}_1 \rightarrow C_0 $, then $$ (\pi_1^\sharp)^{-1} ( \bigcup_j \gamma_j' ) $$ is again a  union 
of paths resulting from disjoint trajectories.  
Now, we iterate this procedure a finite number of times, 
stopping the iteration when the disjoint trajectories in the lift intersect each slit only once.

After each reflection, we restrict to a new minimal width strip, say $C_i$.  Let us call $\tilde{C}_i$ the strip resulting from the $i$th reflection.
We define $\pi^{\sharp}_i :  \tilde{C}_i \rightarrow C_{i-1}$ for $1 \leq i < N$ such that
$$\pi^{\sharp}_i (x,y) = \left\{
 \begin{array}{ll}
 (x,y) & (x,y) \in C_{i-1}, \\
  ( 2 R_i -x,y) & (x,y) \in C_{i-1} '.
  \end{array} \right.$$
Here, $C_{i-1}'$ is defined as the reflected strip and $x=  R_{i-1}$ is the line of reflection for $\tilde{C}_i$.  We can subsequently
define $f_n^{(i)}$ as a sum of delta functions resulting from jumps that occur after reflection, similar to $f_n^{(1)}$ above.
We also have $\pi^N : \R^2 \rightarrow C_N$, the obvious projection that results after we tile the plane with copies of $C_N$.  
So, we have:
$$\R^2 \stackrel{\pi_N}{\rightarrow} C_N \subset \tilde{C}_N \stackrel{\pi^{\sharp}_{N}}{\rightarrow} C_{N-1} \subset \tilde{C}_{N-1}
\stackrel{\pi^{\sharp}_{N-1}}{\rightarrow} ... 
\stackrel{\pi^{\sharp}_{2}}{\rightarrow} C_1 \subset \tilde{C}_1 \stackrel{\pi^{\sharp}_{1}}{\rightarrow} C.$$
Note that 
\[ \pi_N^{-1} (\gamma_1') = \displaystyle\bigcup_j \gamma_{1,j}', \]
where $\{ \gamma_{1,j}' \}$ is the set of all paths in $C_N$ generated by the trajectory $\gamma_1$ and the periodicity in $y$.

\begin{figure}
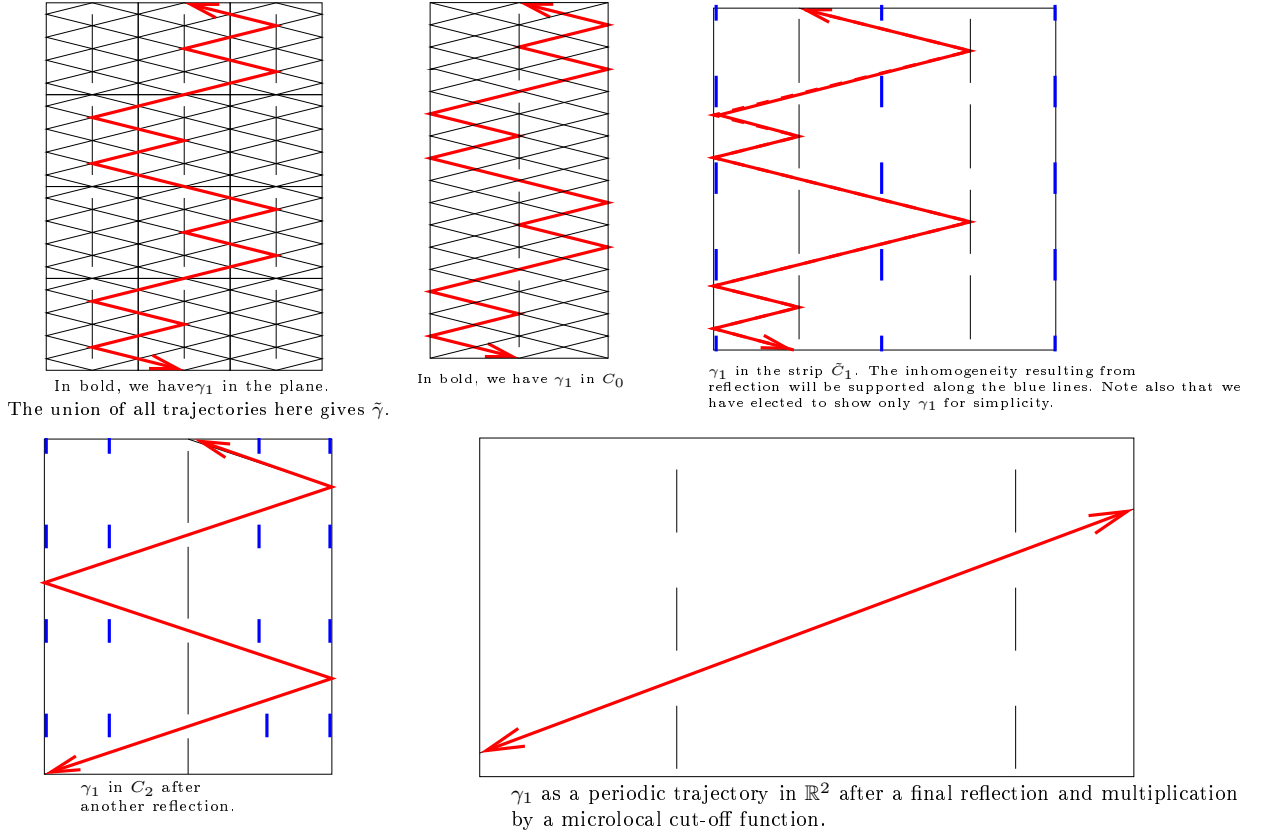

\include{figure6}
\caption{This diagram describes how we "unfold" the eigenfunctions in order to derive a contradiction.}
\label{fig:unf}
\end{figure} 

After a finite number of reflections, we "unfolded" ${\gamma}_1 '$ to be a periodic line on a large strip, $C_N$, which does not 
intersect a slit anywhere.  Now, let us choose $\Phi_\xi$, $\chi$, and $\chi_0$ as above in order to cut-off microlocally 
on this strip around $\gamma_1$.  Again, recall that we can set $\chi \chi_0 = \chi$.  As $f_n^{(i)}$ is supported only in between the slits
for each $i \in \mathbb{N}$, $1 \leq i \leq N$, by
choosing $\Phi_\xi$ to commute with the periodic Laplacian, we have
\[  ( -h^2 \Delta_R - E(h) ) \Phi_\xi \chi {\chi}_0 u_n = \Phi_\xi \chi f_n + [-h^2 \Delta_R , \Phi_\xi \chi ]  u 
=\Phi_\xi  [ -h^2 \Delta_R , \chi ] {\chi}_0 u + {\mathcal O} ( h^{ \infty } ) 
\,, \ \  \| u \| = 1 \,. 
\]
Thus, the result follows by contradiction from the proof of Theorem \ref{t:s}.  
\end{proof}

\begin{rem}
Though this result only shows non-concentration, the proof of Theorem \ref{t:s} can be used to show if $\gamma$ is an $x$-bounded trajectory and
$u$ is an eigenfunction supported on $\gamma ' = \pi_1 (\gamma)$, then in fact there must be mass at the edges of the slits as desired.
\end{rem}

\begin{rem}
If instead of a torus, we had Dirichlet boundary conditions on the boundary of the rectangle as well as the slit, then 
this non-concentration result can also be applied by an elementary reflection principle argument.  
\end{rem}

\end{document}